\documentclass[12pt,reqno,a4paper]{amsart}
\usepackage{graphicx}
\usepackage{amssymb,amscd,textcomp}
\usepackage{hyperref}  

\begin{document}

\let\kappa=\varkappa
\let\epsilon=\varepsilon
\let\phi=\varphi
\let\p\partial
\let\lle=\preccurlyeq
\let\ulle=\curlyeqprec

\def\Z{\mathbb Z}
\def\R{\mathbb R}
\def\N{\mathbb N}
\def\C{\mathbb C}
\def\Q{\mathbb Q}
\def\P{\mathbb P}
\def\HH{\mathsf{H}}
\def\XX{\mathcal X}

\def\conj{\overline}
\def\Beta{\mathrm{B}}
\def\const{\mathrm{const}}
\def\ov{\overline}
\def\wt{\widetilde}
\def\wh{\widehat}

\renewcommand{\Im}{\mathop{\mathrm{Im}}\nolimits}
\renewcommand{\Re}{\mathop{\mathrm{Re}}\nolimits}
\newcommand{\codim}{\mathop{\mathrm{codim}}\nolimits}
\newcommand{\Aut}{\mathop{\mathrm{Aut}}\nolimits}
\newcommand{\lk}{\mathop{\mathrm{lk}}\nolimits}
\newcommand{\sign}{\mathop{\mathrm{sign}}\nolimits}
\newcommand{\rk}{\mathop{\mathrm{rk}}\nolimits}

\def\id{\mathrm{id}}
\def\Leg{\mathrm{Leg}}
\def\Jet{{\mathcal J}}
\def\sS{{\mathcal S}}
\def\lcan{\lambda_{\mathrm{can}}}
\def\ocan{\omega_{\mathrm{can}}}

\renewcommand{\mod}{\mathrel{\mathrm{mod}}}

\newtheorem{mainthm}{Theorem}
\newtheorem{maincor}[mainthm]{Corollary}

\newtheorem{thm}{Theorem}[section]
\newtheorem{lem}[thm]{Lemma}
\newtheorem{prop}[thm]{Proposition}
\newtheorem{cor}[thm]{Corollary}

\theoremstyle{definition}
\newtheorem{exm}[thm]{Example}
\newtheorem{rem}[thm]{Remark}
\newtheorem{df}[thm]{Definition}

\numberwithin{equation}{section}

\title{Instability of rational and polynomial convexity}
\author[Stefan Nemirovski]{Stefan Nemirovski}
\thanks{Partially supported by SFB/TRR 191 `Symplectic Structures in Geometry, 
Algebra and Dynamics', funded by the DFG (Projektnr.\ 281071066~-- TRR 191) 
and by RFBR grant \textnumero 17-01-00592-a}
\address{%
Steklov Mathematical Institute, Gubkina 8, 119991 Moscow, Russia;\hfill\break
\strut\hspace{8 true pt} Mathematisches Institut, Ruhr-Universit\"at Bochum, 44780 Bochum, Germany}
\email{stefan@mi-ras.ru}

\begin{abstract}
It is shown that rational and polynomial convexity of totally real submanifolds
is in general unstable under perturbations that are $C^\alpha$-small for any
H\"older exponent $\alpha<1$. This complements the result of L{\o}w and Wold that 
these properties are $C^1$-stable. 
\end{abstract}

\maketitle

\section{Introduction and results}

A smooth submanifold in $\C^n$ is called totally real if it is not tangent
to any complex line. A compact totally real submanifold is polynomially
or, respectively, rationally convex if and only if every continuous function
on it can be uniformly approximated by polynomials or, respectively, rational 
functions (see e.g.~\cite[Theorem 1.2.10 and Corollary 6.3.3]{Sto}). 

L{\o}w and Wold~\cite[Corollary 1]{LFW} proved that a $C^1$-small perturbation of 
a polynomially convex totally real submanifold $M\Subset\C^n$ remains polynomially
convex. Their argument applies {\it mutatis mutandis\/} to prove that the same
holds for rationally convex totally real submanifolds. The purpose of this note
is to show that these conclusions are no longer true if the perturbations
are only $C^\alpha$-small for any H\"older exponent $\alpha<1$ (and the
dimension of $M$ is $\ge 2$, see Remark~\ref{dim1}).

\begin{df}
\label{perturb-def}
A smooth submanifold $M\subset\C^n$ admits a $C^\alpha$-small
perturbation with a certain property (e.g.\ being rationally convex or not)
if there exist
\begin{itemize}
\item[(i)] a submersion $\pi:U_M\to M$ from a tubular neighbourhood of $M$ onto $M$;
\item[(ii)] a sequence of smooth submanifolds $M_j\subset U_M$, $j\in\N$, having that property 
\end{itemize}
such that $M_j$ converge to $M$ in $C^\alpha$-topology
as $j\to\infty$ and $\pi:M_j\to M$ is a diffeomorphism for all~$j$.
\end{df}

In other words, an identification of a neighbourhood of $M$ in $\C^n$ with a neighbourhood 
of the zero section in its normal bundle can be chosen so that all $M_j$'s correspond 
to graphs of $C^\alpha$-small sections of that bundle. 
For $\alpha\ge 1$ (e.g.\ for deformation considered by L{\o}w and Wold), 
this holds automatically for {\it any\/} such identification once $j$ is large enough,
but in the properly H\"older range this is an additional assumption, cf. Remark~\ref{spintori}(i).

\begin{thm}
\label{main}
Let $T$ be a rationally convex totally real embedded $2$-torus in $\C^2$.
For every $\alpha<1$, there exists a $C^\alpha$-small perturbation of $T$
that is totally real and not rationally convex.
\end{thm}

The tori constructed in the proof of the theorem in~\S\ref{proofthm} 
will be shown not to be rationally convex by an application of a result 
from symplectic geometry on the Maslov class rigidity of Lagrangian
tori~\cite{Pol}. In particular, the perturbed tori are {\it not\/} isotopic 
(or even regularly homotopic) to $T$ through {\it totally real\/} tori
because their Maslov class is different.

\begin{rem}
\label{spintori}
\textbf{(i)} Non rationally convex tori $C^\alpha$-close to a rationally convex one 
can be obtained using {\it spin tori}~\cite[\S 6.4]{DS}. The spin tori 
associated to the embedded curves in $\C\times (0,+\infty)$ 
shown in Fig.~\ref{spin} define a {\it totally real\/} $C^\alpha$-small isotopy 
of the standard product torus in $\C^2$ to a spin torus that is manifestly not rationally convex. 
(Each self-intersection of the rightmost curve in Fig.~\ref{spin} 
corresponds to a holomorphic annulus which has the same oriented
boundary as an annulus on the spin torus and hence lies in its
rationally convex hull by the argument principle.) On the other hand, 
this construction does not seem to yield $C^\alpha$-small perturbations in the sense 
of Definition~\ref{perturb-def} because any fixed normal projection
onto the standard torus will no longer be diffeomorphic on the `kinks' 
once they are small enough.
\begin{figure}[htbp]
\begin{center}
\includegraphics[scale=0.7]{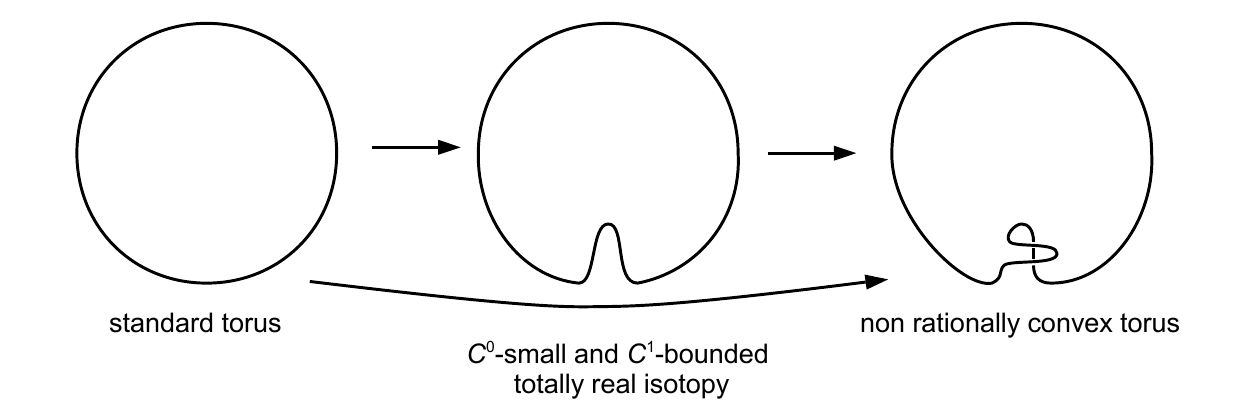}
\end{center}
\caption{Spin tori} 
\label{spin}
\end{figure}

\noindent
\textbf{(ii)} The deformation described in~(i) can be applied to any other totally real surface $S\subset\C^2$
in the following way. First, smoothly approximate $S$ by a real analytic surface~$S'$.
For any embedded real analytic circle $\gamma\subset S'$, there exists a biholomorphic map 
defined on a neighbourhood $U\supset\gamma$ mapping $S'\cap U$ into the product
torus in $\C^2$ and taking $\gamma$ to the fibre of the product torus. Hence, 
the part of the isotopy from~(i) near the `kinks' may be grafted onto~$S'$
making that surface non rationally convex. 
The perturbations in the proof of Theorem~\ref{main} will also be localised 
near an embedded circle on $T$. However, there is no reason to expect 
that grafting them onto other surfaces will have the same effect
because the relevant obstruction to rational convexity is global.
\end{rem}

\begin{cor}
\label{mainc}
The $2$-torus $\wt T=\{(z_1,z_2,z_3)\in\C^3 \mid |z_1|=|z_2|=1, z_1z_2z_3=1\}$
is totally real and polynomially convex in $\C^3$. For every $\alpha<1$, 
there exists a $C^\alpha$-small perturbation of $\wt T$
that is totally real and neither polynomially nor rationally convex.
\end{cor}

\begin{proof} Every continuous function on $\wt T$ can be uniformly approximated
by a Laurent polynomial in $z_1$ and $z_2$ which can be written as a polynomial
in $z_1$, $z_2$, and $z_3=\frac{1}{z_1z_2}$. So $\wt T$ is polynomially convex in $\C^3$.

Let $T^{\mathrm{st}}_j\subset\C^2$, $j\in\N$, be a $C^\alpha$-small perturbation of
the standard torus 
$$
T^{\mathrm{st}}=\{(z_1,z_2)\in\C^2 \mid |z_1|=|z_2|=1\}
$$
that is not rationally convex. Such a perturbation exists by Theorem~\ref{main}. 
Then 
$$
\wt T_j=\{(z_1,z_2,z_3)\in\C^3 \mid (z_1,z_2)\in T^{\mathrm{st}}_j, z_1z_2z_3=1\}
$$
is a $C^\alpha$-small totally real perturbation of~$\wt T$.
If $\wt T_j$ were rationally convex, rational functions of $(z_1,z_2)$
would be dense in the space of continuous functions on~$T^{\mathrm{st}}_j$,
which is impossible since those tori are not rationally convex.
\end{proof}

\begin{rem}
The tori in the corollary {\it are\/} isotopic to $\wt T$ through totally real $2$-tori in $\C^3$. 
This follows by a general position argument. Indeed, the space of complex lines has (real) 
codimension~$\ge 4$ in the space of all real two-planes in $\C^n$ for $n\ge 3$
and therefore a generic isotopy of embedded surfaces in $\C^3$ is through totally
real surfaces.
\end{rem}

Similar examples in higher dimensions are obtained by considering product tori of the form
$$
T_j=\{ z\in\C^n\mid (z_1,z_2)\in T^{\mathrm{st}}_j\text{ and } |z_k|=1 \text{ for }k\ge 3\}
$$ 
and setting
$$
\wt T_j=\{ z\in\C^{n+1}\mid (z_1,\dots,z_n)\in T_j, z_1\cdots z_{n+1}=1\}.
$$
Then $T_j$ is a non rationally convex perturbation of the standard $n$-torus in $\C^n$
and $\wt T_j$ is a non rationally convex perturbation of its polynomially convex lift to $\C^{n+1}$.

\begin{rem}
\label{dim1} 
Polynomial convexity is well-known to be $C^0$-stable for $1$-dimensional totally real tori,
that is, for ($C^1$-)smooth closed curves~$\Gamma\in\C^n$. Indeed, if $\Gamma$ is polynomially 
convex, then every continuous function on $\Gamma$ can be approximated by polynomials and,
in particular, there exist polynomials $P_+$ and $P_-$
such that $\Gamma\cap\{P_+P_-=0\}=\varnothing$ and $\Delta_\Gamma\mathop{\mathrm{arg}} P_\pm =\pm 2\pi$
for some choice of orientation on $\Gamma$. 
A curve $\Gamma'$ sufficiently $C^0$-close to $\Gamma$ is homotopic to it in $\C^n\setminus\{P_+P_-=0\}$
and hence $\Delta_{\Gamma'}\mathop{\mathrm{arg}} P_\pm =\pm 2\pi$  for the induced orientation on $\Gamma'$. 
However, if $\Gamma'$ is not polynomially convex, it is the boundary of a 1-dimensional complex analytic set
by the classical result of Wermer~\cite{Wer}. This implies that for all entire functions the variation 
of the argument along $\Gamma'$ must have the same sign by the argument principle, a contradiction.
\end{rem} 

\begin{rem}
The result of Corollary~\ref{mainc} may be interpreted in the following way. Consider the 2-torus
as the product $S^1\times S^1$ with the angular coordinates $(\theta_1,\theta_2)$.
Laurent polynomials of $e^{i\theta_1}$ and $e^{i\theta_2}$
are dense in the space of continuous functions on the torus. 
The theorem of L{\o}w and Wold asserts that this property
holds for any pair of functions (depending on both variables) that are sufficiently $C^1$-close
to the exponentials. On the other hand, this may no longer be so
for functions that are $C^\alpha$-close to them even if they separate
points on the torus and have $\C$-linearly independent differentials. 
\end{rem}

\section{Proof of Theorem~\ref{main}}
\label{proofthm}
By a result of Duval~\cite{Duv} (see also~\cite{DS}), a totally real surface $\Sigma\subset\C^n$
is rationally convex if and only if it is Lagrangian with respect to a K\"ahler
form $\omega$ on $\C^n$ which can be assumed standard (i.e.~equal to $\omega_{\mathrm{st}}=dd^c\|z\|^2$) 
outside of a sufficiently large ball. (We are only going to use the easier `only if' part of this statement.) 
Let $\omega_t=(1-t)\omega + t\omega_{\mathrm{st}}$, $t\in [0,1]$, be the linear homotopy of K\"ahler forms
connecting $\omega$ to $\omega_{\mathrm{st}}$.
By Moser's stability theorem, there exists a family of compactly supported diffeomorphisms $f_t$, $t\in [0,1]$, of $\C^n$ 
with $f_0=\mathrm{id}$ such that $f_t^*\omega_t=\omega$. Therefore $f_t(\Sigma)$
is a family of $\omega_t$-Lagrangian (and hence totally real) surfaces
connecting $\Sigma=f_0(\Sigma)$ to $f_1(\Sigma)$ which is Lagrangian with respect to $\omega_{\mathrm{st}}$.
It follows, in particular, that if $\Sigma$ is a rationally convex totally real $2$-torus in $\C^2$, 
then its Maslov class $\mu\in H^1(\Sigma;\Z)$ must be non-zero by~\cite[Example 2.2]{Pol}.

Thus, to prove Theorem~\ref{main} it is enough to produce totally real $C^\alpha$-small 
perturbations of $T$ with trivial Maslov class. For $\alpha=0$, the existence 
of such tori may be deduced from Gromov's $h$-principle for totally real embeddings 
(see e.g.~\cite[\S 19.3]{EM}). It is quite probably possible to get H\"older estimates 
from the convex integration proof of the $h$-principle in~\cite{EM}. In the case at hand, 
however, we will use a different approach going back to~\cite{Fie} and based on the cancellation/creation 
theorem for pairs of complex points due to Eliashberg and Kharlamov~\cite{EK}. 

According to~\cite[Example 2.2]{Pol}, the minimal positive value of the Maslov class 
of a Lagrangian torus in $\mathbb C^2$ is always~$2$. Hence, we can find two oriented
simple closed curves $\gamma_1$ and $\gamma_2$ on $T$ intersecting at a single point 
and such that $\mu([\gamma_1])=0$ and $\mu([\gamma_2])=2$. Let $I\subsetneq\gamma_1$ 
be an arc on $\gamma_1$ containing its intersection point with~$\gamma_2$.

\begin{exm}
If $T=T^{\mathrm{st}}$, one can take $\gamma_1=T^{\mathrm{st}}\cap \{z_2=\conj{z}_1\}$
and $\gamma_2=T^{\mathrm{st}}\cap \{z_2=1\}$.
\end{exm}

\begin{figure}[htbp]
\begin{center}
\includegraphics[scale=0.6]{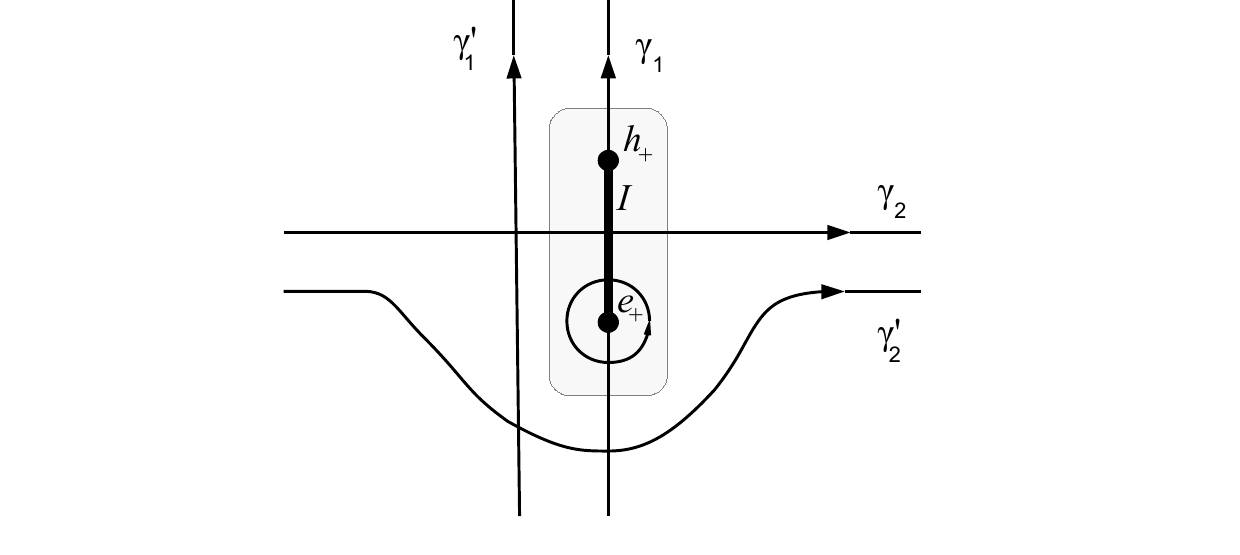}
\end{center}
\caption{First step in the modification of $T$.} 
\label{curves}
\end{figure}

First, we can deform $T$ in a neighbourhood of $I$ so that the new surface has two positive 
(with respect to some orientation on $T$) complex points, one elliptic and one hyperbolic, at the ends of~$I$, see Fig.~\ref{curves}. 
If the elliptic point is on the negative side of $\gamma_2$ with respect to the orientation 
on the torus, then $\mu([\gamma_2])=0$ on the modified surface. (To see this, note that the Maslov class 
takes value $2$ on a small positive loop about the positive elliptic complex point and that this loop gives the homological difference 
in the complement of that complex point between the modified $\gamma_2$ and a representative $\gamma_2'$ of $[\gamma_2]$ on $T$ 
avoiding the neighbourhood affected by the isotopy.) 
Next, we can cancel these two complex points out again by an isotopic modification in a neighbourhood 
of the complementary arc~$\gamma_1\setminus I$. 
The Maslov class of the resulting totally real torus is trivial, as it takes value zero on the 
curves $\gamma_1'$ and $\gamma_2$ generating the first homology group.

To ensure that the modification in the preceding paragraph can be realised by a $C^\alpha$-small perturbation,
we use the argument from the proof of the Eliashberg--Kharlamov theorem given in~\cite[\S 2.3]{Nem} 
and~\cite[\S 10.4]{For}. 
A~neighbourhood $U$ of $\gamma_1$ can be biholomorphically mapped onto a domain in $\C^2$ 
so that $T\cap U$ is taken to the graph of a function $f:A\to\C$ on an annulus~$A\subset\C$. 
(This follows by taking a diffeomorphism of $T\cap U$ onto a model totally real annulus
and approximating it by a holomorphic map.) 
The required modifications are realised by adding to~$f$ 
a suitably cut off solution of the equation $\frac{\p u}{\p\conj z}=\phi$ 
with the right hand side having arbitrarily small $L^p$-norm for any fixed $p>0$. 
(More precisely, $\phi$ may be chosen uniformly bounded by a constant depending 
on the geometry and having arbitrarily thin support.) 
By a classical result (see e.g.~\cite[Theorem 4.3.13]{AIM}),
the Cauchy transform of $\phi$ is a solution of this $\conj\p$-equation 
and its $C^{1-\frac{2}{p}}$-norm is bounded by a constant times the $L^p$-norm of~$\phi$. 
(In~\cite[p.~735]{Nem} and \cite[Lemma 10.4.3]{For},
the $C^0$-version of the same estimate was used.) Thus, $u$ can be made $C^\alpha$-small for any fixed~$\alpha<1$ 
by taking $p$ large enough. 

It should now be clear from the construction that all our deformations are graphs of sections
with respect to a common bundle structure on a neighbourhood of~$T$.\qed

\begin{rem}
It was shown in~\cite{DRGI} that all Lagrangian 2-tori in $\C^2$ with the standard K\"ahler
form are Lagrangian isotopic. So the space of rationally convex totally real 2-tori in $\C^2$ 
is in fact connected. 
\end{rem}

\subsection*{Acknowledgement}
This note is based on the author's talk at the Fourth Central European
Complex Analysis Meeting (Brno, November 2019).

\end{document}